\input epsf.sty
\documentstyle{amsart} 
\newcommand{\al}{\alpha}
\newcommand{\be}{\begin{enumerate}}
\newcommand{\br}{{\Bbb R}}
\newcommand{\bz}{{\Bbb Z}}

\newcommand{\cs}{{\cal S}}
\newcommand{\da}{\delta}
\newcommand{\ee}{\end{enumerate}}

\newcommand{\hyp}{\text{Hyp}^n}
\newcommand{\hyph}{\widehat{\hyp_\la}}
\newcommand{\hypl}{\hyp_\la}

\newcommand{\la}{\lambda}

\newcommand{\lra}{\longrightarrow}
\newcommand{\lw}{\overleftarrow{W$\,\,$}}

\newcommand{\nin}{\noindent}
\newcommand{\pr}{\nin{\bf Proof. }}
\newcommand{\ra}{\rightarrow}
\newcommand{\res}{\text{Res}\,}

\newcommand{\rw}{\overrightarrow{W$\,\,$}}

\newcommand{\sm}{\setminus}

\newcommand{\susp}{\text{susp}\,}

\newcommand{\ti}{\tilde}
\newcommand{\wti}{\widetilde}

\newtheorem{thm}{Theorem}[section]
\newtheorem{sthm}[thm]{First Stabilization Theorem}
\newtheorem{sthm2}[thm]{Second Stabilization Theorem}
\newtheorem{df}  [thm]{Definition}

\newtheorem{crl} [thm]{Corollary}
\newtheorem{prop}[thm]{Proposition}
\newtheorem{conj}[thm]{Conjecture}

\numberwithin{equation}{section}

\begin{document}

\title[Topology of spaces of hyperbolic polynomials]
{Topology of spaces of hyperbolic polynomials and combinatorics of resonances}
              \author{Dmitry N. Kozlov}
               \date{\today}

   \address{{\it current address:} Department of Mathematics,
University of Bern, CH-3012, Bern, Switzerland}
\address {{\it on leave from:} Department of Mathematics, 
Royal Institute of Technology,
S-100 44, Stockholm, Sweden}
\email{ kozlov@@math.ethz.ch,kozlov@@math.kth.se}

\begin{abstract} 
  In this paper we study the topology of the strata, indexed by number 
partitions $\lambda$, in the natural stratification of the space of monic
hyperbolic polynomials of degree~$n$. We prove stabilization theorems
for removing an independent block or an independent relation in $\lambda$.
We also prove contractibility of the one-point compactifications
of the strata indexed by a~large class of number partitions, 
including $\lambda=(k^m,1^r)$, for $m\geq 2$. Furthermore, we study the
maps between the homology groups of the strata, induced by imposing
additional relations (resonances) on the number partition $\lambda$, 
or by merging some of the blocks of $\lambda$. 
\end{abstract}

\maketitle

                 \section{Introduction}

   The space of monic hyperbolic polynomials of degree $n$ is naturally
stratified by fixing the multiplicities of the roots. In this paper 
we study the topology of these strata. The topological aspects of the 
spaces of polynomials with multiple roots have been extensively studied,
see e.q., \cite{Ar1,Ar2,SW,V1,V2,V3}. As the general motivation we would 
like to mention the widely-branching program of studying the topological
properties of certain subsets of some fixed function space, namely of
the spaces of functions with singularities of some fixed type 
(also known as {\it discriminants}). 

  Here we study the discriminants in the function spaces of all
hyperbolic polynomial maps from $\br^n$ to $\br^n$, but there is 
a very large number of other important examples:
\begin{itemize}
\item spaces of knots, i.e., nonsingular imbeddings $S^1\lra S^3$, e.g., 
see~\cite{V4};
\item spaces of complex polynomial maps, or, more generally, of systems
of polynomials, see for example \cite[Theorem 4, p.~126]{V1}, and \cite{CCMM}, 
for a~connection between these two cases;
\item spaces of nonsingular deformations of differentiable manifolds; 
\end{itemize}
 and many others. The V.~Vassiliev's book $\cite{V1}$ provides a~very 
extensive overview of the developments and the current state of the art
in this area.

  The case of hyperbolic polynomials has been considered before, most
importantly in \cite{SW}, where a~simplicial complex of a~combinatorial nature
$\da_\la$ was described, such that the double suspension of $\da_\la$ is
homeomorphic to $\hyph$. Here $\hyph$ is a~one-point compactification of 
the strata of the space of all monic hyperbolic polynomials of degree $n$,
which is indexed by the number partition $\la$; the exact definition is 
given in Subsection~\ref{ss2.1}.

  In this paper we are working further with this combinatorial model. 
By using the techniques of the Discrete Morse Theory as well as direct
algebro-topological arguments, we are able to compute homology groups
of $\hyph$ and even, in some cases, determine the homotopy type of~$\hyph$,
for several previously unknown classes of~$\la$. 

We would like to emphasize the combinatorial aspect at this point. 
The methods which we use are combinatorial. In fact, we are working
exclusively with posets of compositions, which can be considered as
objects of internal combinatorial interest, even though their
appearance was mainly motivated by the topological questions about
$\hyph$.

Boris Shapiro has suggested to me in private conversation, \cite{Sh}, that 
there might exist a~general algorithm for computing the homology groups (or 
even better, the homotopy type) of $\hyph$ for general $\la$. The results of 
this paper may be used as the first step on this path.

\vskip3pt
  Here is a~short summary of the contents.

\vskip3pt
\noindent
{\bf Section 2}. Notations and terminology are introduced. The
description of the Shapiro-Welker combinatorial model for $\hyph$ is
given. The relevant results of the Discrete Morse Theory are outlined
for later use.

\vskip3pt
\noindent
{\bf Section 3}. Here the bulk of our results is concentrated.  In
Subsection~\ref{ss3.1}, we prove the First Stabilization Theorem which
allows one to remove an independent block from a~number partition
$\la$.  Furthermore, the Theorem~\ref{chop} describes, in particular,
a~large class of partitions $\la$ (to which for example
$\la=(k^m,1^r)$, $m\geq 2$, belongs), for which $\hyph$ is
contractible.

In Subsection~\ref{ss3.2}, we prove the Second Stabilization Theorem
which allows one in some situations to remove a~relation from
a~partition. We also study the maps between simplicial complexes
$\da_\la$ which are induced by imposing additional relations
(resonances) on the number partition~$\la$, or by merging some of the
blocks of $\la$.

\vskip3pt
\noindent
{\bf Section 4}. We describe some remaining open questions and perspective
developments.

\vskip4pt
\noindent
{\bf Acknowledgments.} I am grateful to Boris Shapiro for educational
discussions and for inspiring this research, and to Eva-Maria
Feichtner for the careful proofreading of the early versions of the
paper. I would also like to thank the referee for helpful comments
which led to substantial improvements throughout the paper.
Finally, I acknowledge the support of the Forschungsinstitut f\"ur
Mathematik, ETH Z\"urich, where this paper was written.

         \section{Methods}
   
\subsection{Compositions, number partitions, and the indexing of strata} 
\label{ss2.1} $\,$
\vskip3pt

An ordered tuple of positive integers $\la=[\pi_1,\dots,\pi_t]$ is
called a~{\it composition}, or, sometimes a~{\it composition of $n$},
where $n=\pi_1+\dots+\pi_t$. When this tuple is taken unordered, $\la$
is called a~{\it number partition} of $n$, we write,
$\la=(\pi_1,\dots,\pi_t)$, and $\la\vdash n$. For number partitions,
we also use the power notation: $(n^{\al_n},\dots,1^{\al_1})=
(\underbrace{n,\dots,n}_{\al_n},\dots,\underbrace{1,\dots,1}_{\al_1})$.
Both for number partitions and compositions, we call $\pi_i$'s the
{\it blocks} of~$\la$. The {\it length} of $\la$ is the number of
blocks, it is denoted~$l(\la)$. Given a~composition of $n$, its {\it
  type} is the number partition of $n$, which is obtained from the
composition by forgetting the order of the blocks (in the text we
often reflect it by changing the square brackets to the round ones).

Let $\hyp\subseteq\br^n$ be the space of all monic hyperbolic
polynomials (a polynomial is called hyperbolic when all of its roots,
and hence its coefficients, are real numbers), and $\widehat\hyp$ its
one-point compactification. There is a~standard cell decomposition of
$\widehat\hyp$ which we now proceed to describe. For a~composition
$[\al_1,\dots,\al_t]$ of $n$, we denote by
$\hyp_{[\al_1,\dots,\al_t]}$ the topological space of all hyperbolic
polynomials $(x-r_1)^{\al_1}\dots(x-r_t)^{\al_t}$ such that
$r_1<\dots<r_t$. Given a~number partition $\la$ of $n$, we denote by
$\hypl$ the closure (in $\hyp$) of the union of all cells
$\hyp_{[\al_1,\dots,\al_t]}$, where the composition
$[\al_1,\dots,\al_t]$ is of type~$\la$. We denote the one-point
compactification of $\hypl$ by $\hyph$.

 \subsection{Shapiro-Welker model}$\,$
\vskip3pt

The set of all compositions of $n$ is partially ordered by refinement.
Namely, let $x=[\al_1,\dots,\al_{l(x)}]$ and $y=[\beta_1,\dots,\beta_{l(y)}]$
be two compositions of $n$, we say that $x\leq y$ if and only if
$\al_j=\beta_{i_{j-1}+1}+\dots+\beta_{i_j}$, for $1\leq j\leq l(x)$, and
some $0=i_0<i_1<\dots<i_{l(x)}=l(y)$. Since $\beta_i>0$, for $i=1,\dots,l(y)$,
the indices $i_1,\dots,i_{l(x)-1}$ are uniquely defined. In this situation, 
we set $g(y,x,\beta_i)=j$ if and only if $i_{j-1}+1\leq i\leq i_j$.

  Given a~number partition $\la=(\pi_1,\dots,\pi_t)$ of $n$, we define
$D_\la$ to be the poset consisting of all compositions of $n$ which 
are less or equal of some composition of $n$ of type $\la$. $D_\la$
has a~minimal element, the composition consisting of just the number $n$,
and it is easy to see that $D_\la\cup\{\hat 1\}$ is a~lattice, where 
$\hat 1$ is an adjoint maximal element.  

  Since the lower intervals of $D_\la$ are boolean algebras, and $D_\la$ 
itself is meet-semilattice, there exists a~unique simplicial complex, which 
we denote by $\da_\la$, such that $D_\la$ is the face poset of $\da_\la$,
i.e., the elements of $D_\la$ and the simplices of $\da_\la$ are in 
bijection, and the partial order relation on $D_\la$ corresponds under
this bijection to the inclusions of simplices of $\da_\la$. In particular, 
when $l(\la)=1$, we have $\da_\la=\emptyset$. 

This bijection is the reason for why we chose an~order convention on
the compositions opposite to the customary in combinatorics: we want
to have an~order-preserving bijection, not an order-reversing one.

The simplicial complex $\da_\la$ is important for the following reason. 

\begin{thm} (\cite[Theorem 3.5(a)]{SW}).
Let $\la$ be a~number partition of $n$, then the one-point compactification
of the strata indexed by $\la$, $\hyph$, is homeomorphic
to the double suspension of the simplicial complex $\da_\la$.
\end{thm}

\subsection{Terminology of resonances}$\,$
\vskip3pt

It turns out that $\da_\la$ depends only on the set of various
equalities of sums of different parts of $\la$ (resonances), not on
the exact numerical values of the parts of~$\lambda$. This can be
formalized as follows. 

\begin{df} $\,$
  
\nin {\bf (a)} Given a~composition of $n$,
$\alpha=[\pi_1,\dots,\pi_{l(\alpha)}]$, a~{\bf resonance} 
of $\alpha$ is an unordered pair $(\{i_1,\dots,i_k\},\{j_1,\dots,j_m\})$
of nonempty disjoint subsets of $\{1,\dots,l(\alpha)\}$, such that
\begin{equation}
  \label{eq:2.1}
  \pi_{i_1}+\dots+\pi_{i_k}=\pi_{j_1}+\dots+\pi_{j_m}.
\end{equation}
 
We denote the set of all resonances of $\alpha$ by $\res\alpha$.

\nin {\bf (b)}
  Given a~number partition $\la\vdash n$, we denote by $\res\la$
the set of all $\res\alpha$, such that $\alpha$ has type $\lambda$.

\end{df}

For every positive integer $k$, the permutation action of $\cs_k$ on
$[k]$ induces an $\cs_k$-action on the set of all compositions of
length $k$, and hence also on the set $\{\res\alpha\,|\,\alpha$ is
a~composition of length~$k \}$. The orbits are indexed by number
partitions, and $\res\la$ is the orbit of this action indexed by
$\lambda$. For example $\res(8,3,3,3,1)=\res(7,5,4,4,4)$.

The notion of a~resonance is important, because if $\lambda$ and
$\tilde\lambda$ have the same set of resonances, then the spaces
$\hyph$ and $\widehat{\text{Hyp}_{\tilde\lambda}^n}$ are homeomorphic.

We will abuse the language and call an~equality of the
form~\eqref{eq:2.1} itself a~resonance. We will also say that $\alpha$
(or $\lambda$) has this resonance.  Where it does not lead to the
confusion, we shall often say "the set of resonances of the number
partition $\lambda$", instead of saying "the set of resonances of some
composition of type $\lambda$". If in~\eqref{eq:2.1} $k=m=1$, then the
resonance is called {\it trivial}, otherwise it is called {\it
  nontrivial}.  The resonance
$\pi_{i_1}+\dots+\pi_{i_k}=\pi_{j_1}+\dots+\pi_{j_m}$ is said {\it to
  involve} blocks $\pi_{i_1},\dots,\pi_{i_k},
\pi_{j_1},\dots,\pi_{j_m}$, correspondingly these blocks are said {\it
  to be involved} in this resonance.

A block $\pi_i$ of a~composition is called {\it independent}, if,
whenever $\pi_{i_1}+\dots+\pi_{i_k}=\pi_{j_1}+\dots+\pi_{j_m}$, and
$i_1=i$, there exists $1\leq q\leq m$, such that there is a~trivial
resonance $\pi_i=\pi_{j_q}$. If a~block is not involved in any
resonance than it is called {\it strongly independent}.

Clearly, given a~block in a~number partition $\la$, the corresponding
block is independent, resp. strongly independent, in all compositions
of type $\lambda$ if and only if it is independent, resp. strongly
independent, in any one such composition. Therefore, we have
a~well-defined notion for a~block of a~number partition to be
independent, resp. strongly independent.

If $\la=(\pi,\pi_1,\dots,\pi_{l(\la)-1})$ is a~number partition such
that $\pi$ is strongly independent, then for any $x\in D_\la$ there
exists a~number $1\leq\rho_\la(x,\pi)\leq l(x)$, such that, whenever
$y\geq x$ and $y$ is of type~$\la$, we have
$\rho_\la(x,\pi)=g(y,x,\pi)$, i.e., this number does not depend on the
choice of $y$.

More generally, if $\pi$ is independent, but not necessarily strongly
independent: say the trivial resonances in which it is involved are
$\pi=\pi_1$, $\dots$, $\pi=\pi_k$, then the multiset
$\{g(y,x,\pi_i)\,|\,i=1,\dots,k\}\cup\{g(y,x,\pi)\}$ does not depend
on the choice of the composition $y$ of type $\la$, such that $y\geq
x$. In such situation, if $x=[\al_1,\dots,\al_{l(x)}]$, we say that
$\al_j$ {\it contains the block} $\pi$, if $j$ belongs to the multiset
$\{g(y,x,\pi_i)\,|\,i=1,\dots,k\}\cup\{g(y,x,\pi)\}$. We refer to the
number of occurrences of $j$ in this multiset as the number of copies
of $\pi$ contained in $\al_j$.

 \subsection{Discrete Morse theory} \nopagebreak[4]$\,$
\nopagebreak[4]\vskip3pt\nopagebreak[4]

For a~regular CW complex $\da$, we denote by $P(\da)$ its face poset,
the empty face included as a~minimal element. Vice versa, for
a~poset~$P$, $\da(P)$ denotes the simplicial complex which is the
nerve (order complex) of $P$, see \cite{Q,Se} for the first
appearances of the nerve of a~category. Recall also the following
terminology.

\begin{df}\label{coll}
  Let $X$ be a~regular CW complex. Assume that $F_1$ and $F_2$ are cells of $X$ 
such that $F_2$ is a~maximal cell which contains $F_1$, and there is no other
ma\-xi\-mal cell containing $F_1$. A~{\bf collapse} is the replacement of $X$ 
with $X\sm\{F\,|\,F_1\subseteq F\}$. A~collapse is called {\bf elementary} if 
$\dim F_1+1=\dim F_2$. 
\end{df}

  Clearly, a~collapse is a~strong deformation retract, hence it preserves the 
homotopy type of the space.

  Let $\da$ be a~regular CW complex. A~{\it matching} $W$ on $P=P(\da)$ 
(cf.~\cite[Definition~9.1]{Fo}) is a~set of disjoint pairs $(\sigma,\tau)$ such 
that $\tau,\sigma\in P$, $\tau\succ\sigma$, (``$\succ$'' denotes the covering 
relation). We set 
$$\rw=\{\sigma\in P\,|\, \text{ there exists } \tau \text{ such that } 
(\sigma,\tau)\in W\},$$ 
$$\lw=\{\tau\in P\,|\, \text{ there exists } \sigma \text{ such that } 
(\sigma,\tau)\in W\}.$$ 
If $(\sigma,\tau)\in W$ then we set $W(\sigma)=\tau$.

\begin{df} (cf. \cite[Definition~9.2]{Fo}).
  A~matching is called {\bf acyclic} if it is impossible to find a~sequence 
$\sigma_0,\dots,\sigma_t\in\rw$ such that $\sigma_0\neq \sigma_1$, 
$\sigma_0=\sigma_t$ and $W(\sigma_i)\succ\sigma_{i+1}$ for $0\leq i\leq t-1$.
\end{df}
Note, that if a~matching is acyclic, then not all 0-cells are matched
with 1-cells. 

A cell $\sigma$ is called {\it critical} if, either it is
the empty cell, or it is a~0-cell matched with the empty cell, or
$\sigma\not\in\lw\cup\rw$; in the latter case $\sigma$ is called
{\it nontrivial critical}. Let $m_i(W)$ denote the number of critical
$i$-cells.

\vskip3pt
\nin{\bf Note.} Alternatively, we could just omit the empty cell from
the cell complex. However, we choose to keep it, and have a~somewhat
more complicated definition of the critical cells. First, since the
empty cell is natural in our applications, and, second, since the 
formulation of the Theorem~\ref{morse} is somewhat smoother in that
version, as a~complete matching is a~standard object in combinatorics. 

\vskip3pt
 We need the following result, see also \cite[Theorem 3.2]{Ko}, 
\cite[Theorem 9.3]{Fo}, and \cite[Corollary~3.5]{Fo}.

\begin{thm}\label{morse}
  Let $\da$ be a~regular CW complex of dimension $d$, and let $W$ be
  an acyclic matching on $P(\da)$. Then $\da$ is homotopy equivalent
  to a~CW complex $\da^M$, which has $m_i(W)$ cells of dimension $i$.
  In particular, if the acyclic matching is complete, then $\da$ is
  contractible.
\end{thm}

The basic idea of the proof is that the combinatorial condition
of acyclicity allows us to arrange the collapses in a~sequence
and perform them one after the other.

\section{Results}

\subsection{Applications of Discrete Morse theory}$\,$
\label{ss3.1}
\vskip3pt

  We begin by proving a~theorem which allows one to remove an independent
block from a~number partition.

\begin{sthm} \label{stab}
  Let $\la=(\pi_1,\dots,\pi_t)$ be a~number partition of~$n$, 
such that $\pi_1$ is independent.
\begin{enumerate}
\item[(a)] If $\pi_1$ is not strongly independent, i.e., $\pi_1=\pi_i$, 
for some $2\leq i\leq t$, then the simplicial complex $\da_\la$, 
and therefore also the topological space $\hyph$, is contractible. 
\item[(b)] If, on the other hand, $\pi_1\neq \pi_i$, for all $2\leq i\leq t$, 
then $\da_\la$ is homotopy equivalent to $\susp\da_{\ti\pi}$,
correspondingly $\hyph\simeq\susp\widehat{\text{Hyp}_{\ti\pi}^n}$,
where $\ti\pi=(\pi_2,\dots,\pi_t)$. 
\end{enumerate}
\end{sthm}
\noindent \pr

\noindent 
{\bf (a)} For any $a=[\al_1,\dots,\al_{l(a)}]\in D_\la$ let $\iota(a)$, 
resp.~$\gamma(a)$, be the smallest index of a~block of $a$ which is equal 
to $\pi_1$, resp.~larger than $\pi_1$ and containing at least one copy of 
$\pi_1$, if such exists, otherwise put $\iota(a)$, resp. $\gamma(a)$, 
equal to $\infty$. Clearly at least one of the numbers $\iota(a)$ and 
$\gamma(a)$ is finite. Since 
$\iota(a)\neq\gamma(a)$, the elements of $D_\la$ split into two disjoint sets: 
$A=\{a\in D_\la\,|\,\iota(a)>\gamma(a)\}$ and 
$B=\{a\in D_\la\,|\,\iota(a)<\gamma(a)\}$.

Consider the following matching $W$ on $D_\la$. For $a=[\al_1,\dots,\al_{l(a)}],
b=[\beta_1,\dots,\beta_{l(b)}]\in D_\la$, $a\in A$, $b\in B$, we have 
$(a,b)\in W$ if and only if $l(a)+1=l(b)$, and 
$$\begin{cases}
   \al_i=\beta_i, & \text{ for } i=1,\dots,\iota(b)-1;\\
   \al_{\iota(b)}=\beta_{\iota(b)}+\beta_{\iota(b)+1};&\\
   \al_i=\beta_{i+1},& \text{ for } i=\iota(b)+1,\dots,l(a).
\end{cases} $$ 
Note that in such case $\iota(b)=\gamma(a)$, and $\iota(b)\leq l(b)-1$.

  Given $a\in A$, by the definition of W, there exists a~uniquely defined
$b\in B$, such that $(a,b)\in W$; reversely, given $b\in B$, one obtains
a unique $a$ such that $(a,b)\in W$ by breaking the block $\beta_{\gamma(b)}$
into the blocks $\pi_1,\beta_{\gamma(b)}-\pi_1$ (exactly in this order).

It is clear that these procedures are inverses of each other, thus $W$
is a~well-defined matching. Furthermore, $W$ is complete. Note that
the empty simplex $[\pi_1+\dots+\pi_t]$ is matched with the vertex
$[\pi_1,\pi_2+\dots+\pi_t]$. We have $\rw=A$ and $\lw=B$.
 
The topologically inclined reader may think of this matching as
a~set of elementary collapses. We shall now check that they in fact
can be arranged in a~sequence of collapses, by checking that $W$ is an
acyclic matching. Consider a~sequence $a_0,\dots,a_k\in A$, such that
$a_0=a_k$, $a_i\neq a_{i+1}$ and $W(a_i)\succ a_{i+1}$, for
$i=0,\dots,k-1$. The simple but crucial observation is that
$\iota(W(a_i))>\gamma(a_{i+1})$.

Indeed, since $a_{i+1}\in A$ and $W(a_i)\in B$, we have 
$\gamma(a_{i+1})<\iota(a_{i+1})$ and $\gamma(W(a_i))>\iota(W(a_i))$.
Therefore, there are two cases. The first one is when $a_{i+1}$ is obtained 
from $W(a_i)$ by merging blocks with indices $\iota(W(a_i))-1$ and 
$\iota(W(a_i))$, in which case $\gamma(a_{i+1})=\iota(W(a_i))-1$. 
The second case is when $a_{i+1}$ is obtained from $W(a_i)$ by merging 
blocks with indices $j-1$ and $j$, such that $j<\iota(W(a_i))$, and the 
obtained block is larger than $\pi_1$ and contains at least one copy of $\pi_1$, 
in which case $\gamma(a_i)=j-1<\iota(W(a_i))$. Note that we are using the 
fact that $a_i\neq a_{i+1}$ by ruling out the possibility that $a_{i+1}$ 
is obtained from $W(a_i)$ by merging the blocks indexed 
$\iota(W(a_i))$ and $\iota(W(a_i))+1$. 

Since, as observed before, $\gamma(a_i)=\iota(W(a_i))$, we obtain 
$\gamma(a_i)>\gamma(a_{i+1})$ and hence a~contradiction 
$\gamma(a_0)>\gamma(a_1)>\dots>\gamma(a_t)=\gamma(a_0)$.
Thus, by the Theorem \ref{morse}, the simplicial complex $\da_\la$ is contractible.

\noindent
{\bf (b)} Let $x$ be the vertex of $\da_\la$ which is labeled by the
composition $[\pi_2+\dots+\pi_t,\pi_1]$. Consider the two following
subspaces of $\da_\la$. $T=\{{\text{closed star}\,}_{\da_\la}(x)\}$
and $Q=\da_\la\setminus \{{\text{open star}\,}_{\da_\la}(x)\}$.
Clearly $T$ is contractible and $\da_\la=T\cup Q$.

Let us show that $Q$ is contractible. Consider the matching $W$ on
$\da_\la$ which is defined completely analogously to the one in the
first part of this proof. It is easy to see that the matching is again
complete. The acyclicity of $W$ follows from the argument in the first
part.

  Then $\da_\la$ is homotopy equivalent to the suspension of 
$T\cap Q={\text{link}\,}_{\da_\la}\,(x)$. The simple observation that 
${\text{link}\,}_{\da_\la}(x)=\da_{(\pi_2,\dots,\pi_t)}$ finishes the proof. 
\qed

\vskip3pt
  Shapiro and Welker have computed the homotopy type of $\da_\la$ for 
$\la=(k,1^r)$, \cite[Proposition 3.9, Corollary 3.10]{SW}, by using the previous 
work of Bj\"orner and Wachs, \cite[Theorem 8.2, Corollary 8.4]{BW}.
Here we prove a~general theorem of which $\la=(k^m,1^r)$, for $m\geq 2$,  
is a~special case. We shall also reprove the result for the case $\la=(k,1^r)$.

First we need some terminology. Number partitions of $n$ are partially
ordered by refinement: for $\la,\mu\vdash n$,
$\la=(\pi_1,\dots,\pi_{l(\la)})$, $\mu=(\eta_1,\dots,\eta_{l(\mu)})$,
we say that $\la\geq\mu$ if there exists a~collection of disjoint sets
$I_1,\dots,I_{l(\mu)}\subseteq\{1,\dots,l(\la)\}$, such that
$\cup_{i=1}^{l(\mu)}I_i=\{1,\dots,l(\la)\}$ and $\sum_{j\in
  I_k}\pi_j=\eta_k$, for $k=1,\dots,l(\mu)$. Again, the convention of
the ordering is dictated by the partial order on the set of
compositions, which in turn followed the pattern of cell inclusions.

\vskip3pt
 {\bf Cutting Condition.} {\it We say that a~pair $(\la,\pi_1)$, where 
$\la=(\pi_1,\dots,\pi_{l(\la)})$ is a~number partition of $n$, satisfies 
the Cutting Condition, if, whenever $\mu=(\eta_1,\dots,\eta_{l(\mu)})$ 
is another number partition of $n$, such that $\la\geq\mu$, and for some 
$i\in\{1,\dots,l(\mu)\}$ and some nonempty set $I\subseteq\{2,\dots,l(\la)\}$, 
we have equality $\eta_i=\pi_1+\sum_{j\in I}\pi_j$, then we have $\la\geq\ti\mu$, 
where $\ti\mu=(\eta_1,\dots,\eta_{i-1},\pi_1,\sum_{j\in I}\pi_j,\eta_{i+1},
\dots,\eta_{l(\mu)})$.}
\vskip3pt

  Note that if $(\la,\pi_1)$ satisfies the cutting condition, then 
the block $\pi_1$ is not necessarily the largest one. For example, 
$((6,4,4,2,1),1)$ satisfies the cutting condition, while
$((6,4,4,2,1),6)$ does not: $(6,4,4,2,1)>(8,6,2,1)$, but 
$(6,4,4,2,1)\not\geq(6,6,2,2,1)$.

\begin{thm} \label{chop}
  Let $\la=(\pi_1,\dots,\pi_{l(\la)})$ be a~number partition, such that 
$(\la,\pi_1)$ satisfies the cutting condition.
\begin{enumerate}
\item [(a)] if $\pi_1$ is involved in a~trivial resonance, then
  $\da_\la$ is contractible, and hence so is $\hyph$;
\item [(b)] if $\pi_1\neq\pi_i$, for $i=2,\dots,l(\la)$, then there
  exists an acyclic matching on $D_\la$, such that the nontrivial
  critical simplices are exactly the ones indexed by the compositions
  $a=(\al_1,\dots,\al_{l(a)})$, where $\al_{l(a)}=\pi_1$, and
  $\al_i\neq\pi_1+\sum_{j\in I}\pi_j$, for all $1\leq i\leq l(a)-1$
  and $I\subseteq\{2,\dots,l(\la)\}$ ($I$ may be empty).
\end{enumerate}
\end{thm}

\pr We can define a~matching on $D_\la$ which is completely analogous
to the one defined in the proof of the Stabilization
Theorem~\ref{stab}(a). A~word by word check of the proof of the
acyclicity of the matching reveals that the same argument is still
valid in our case.

If $\pi_1\neq\pi_i$, for $i=2,\dots,l(\la)$, the nontrivial critical
simplices are the simplices corresponding to the compositions
described in (b) above, because they are the only ones where on one
hand $\pi_1$ cannot be merged with the next block to the right, and on
the other hand, it is impossible to cut off a~block $\pi_1$ from some
other block. 
\qed

The Theorem~\ref{stab}(a) follows from the Theorem~\ref{chop}(a).
However, since the proofs are essentially identical, we prefer to
prove the structural Stabilization Theorem~\ref{stab} first, and then
point out that the argument is actually valid for the combinatorially
more technical Theorem~\ref{chop}. The next result shows that the
Theorem~\ref{chop} is strictly more general that the
Theorem~\ref{stab}.

\begin{crl}
  Let $\la=(\underbrace{k^{c_1},\dots,k^{c_1}}_{m_1},\dots,
  \underbrace{k^{c_t},\dots,k^{c_t}}_{m_t},\underbrace{1,\dots,1}_{m_{t+1}})$,
  for some positive integers $k,m_1,\dots,m_{t+1}$, and
  $c_1>\dots>c_t$, such that $k\geq 2$, and $m_1\geq 2$, then
  $\da_\la$, and hence also $\hyph$, is contractible.
\end{crl}

\pr It is enough to check that $(\la,k^{c_1})$ satisfies the cutting
condition, since then, by the Theorem~\ref{chop}(a), the simplicial
complex $\da_\la$ is contractible.

  Assume $\mu=(\eta_1,\dots,\eta_{l(\mu)})$, $\la>\mu$, and 
$\eta_q=\sum_{i=1}^t r_i k^{c_i}+r_{t+1}$ for $r_i\leq m_i$, $r_1\geq 1$,
and either $r_1\geq 2$ or $r_i>0$ for some $i\in\{2,\dots,t+1\}$. Since 
$\la\geq\mu$, one can write $\eta_j=\sum_{i=1}^t r_{i,j} k^{c_i}+r_{t+1,j}$, for 
$j=1,\dots,l(\mu)$, so that $\sum_{j=1}^{l(\mu)}r_{i,j}=m_i$, for $i=1,\dots,t+1$. 
If $r_{1,q}\geq 1$, then we are done. Otherwise, as $\eta_q>k^{c_1}$, 
we can find $\ti r_i\leq r_{i,q}$, for $i=2,\dots,t+1$, such that 
$\sum_{i=2}^t \ti r_i k^{c_i}+\ti r_{t+1}=k^{c_1}$. This means that there 
exist numbers $\ti r_{i,j}$ such that 
$\eta_j=\sum_{i=1}^t \ti r_{i,j} k^{c_i}+\ti r_{t+1,j}$, for $j=1,\dots,l(\mu)$,
$\sum_{j=1}^{l(\mu)}\ti r_{i,j}=m_i$, for $i=1,\dots,t+1$,
and $\ti r_{1,q}\geq 1$. 
\qed

\vskip3pt
 There are many other pairs satisfying the cutting condition, for example:
\begin{itemize}
\item $((\pi^p,\pi_1^{q_1},\dots,\pi_m^{q_m},1^r),\pi)$, for 
$\pi\geq\sum_{i=1}^m q_i\pi_i$;
\item $((k\pi_1,\dots,k\pi_t,1^r),1)$, for $1\leq r\leq k-1$.
\end{itemize}

Often the matching produced in the Theorem~\ref{chop}(b) can be 
extended so as to yield the complete information on the homology
groups of the simplicial complex~$\da_\la$. In the next proposition
we demonstrate this on a~couple of examples.

\begin{prop}\nopagebreak[4] \label{pr3.4}\nopagebreak[4]$\,$\nopagebreak[4]
\begin{enumerate}\nopagebreak[4]
\item [(a)] (\cite[Proposition 3.9, Corollary 3.10]{SW}). For $\la=(k,1^t)$, 
where $k\geq 2$, $t\geq 0$, we have 
$$\da_\la\simeq\begin{cases}
  S^{2m-1},& \text{ if } t=km, \text{ for some } m\in{\Bbb Z};\\
  S^{2m},& \text{ if } t=km+1, \text{ for some } m\in{\Bbb Z};\\
  \text{point},& \text{ otherwise}.
\end{cases}$$ 
\item [(b)] For $\la=(k,2,1^t)$, where $k\geq 3$, we have 
$$\wti H_*(\da_\la)=
\begin{cases}
\bz_{(t)}\oplus\bz_{(2m-1)},&\text{ if } t+2=km, \text{ for some } m\in{\Bbb Z};\\
\bz_{(t)}\oplus\bz_{(2m)},&\text{ if } t+1=km, \text{ for some } m\in{\Bbb Z};\\
\bz_{(t)},&\text{ otherwise.}
\end{cases}$$
\end{enumerate}
Recall that $\hyph\cong\susp^2\da_\la$, hence (a) and (b) above
yield the corresponding information about the stratum $\hyph$.
\end{prop}

\nin {\bf Note.} In this paper, all homology groups are reduced and
with integer coefficients. We also use the notation $\bz_{(i)}$ to
denote a~direct summand $\bz$ in the $i$th reduced homology group.
For example, the reduced homology groups of the torus $S^1\times S^1$
would be written as $\bz_{(1)}\oplus\bz_{(1)}\oplus\bz_{(2)}$.

\vskip3pt
\nin{\bf Proof of Proposition~\ref{pr3.4}.}

\vskip3pt
\noindent
{\bf (a)} Let us extend the matching given in the proof of the
Theorem~\ref{chop}(b) (equivalently, in the proof of the Stabilization
Theorem~\ref{stab}) as follows. If $a\in D_\la$ indexes a~nontrivial
critical simplex and $a\neq[1,k-1,1,\dots,k-1,1,k]$,
$a\neq[1,k-1,\dots,1,k-1,k]$, then
$a=[\underbrace{1,k-1,\dots,1,k-1}_{2m},p,q,\dots,k]$, where $m\geq
0$, $p\leq k-1$, and, either $p\geq 2$ or $q\leq k-2$.

If $p\geq 2$, we define $W(a)=[\underbrace{1,k-1,\dots,1,k-1}_{2m},
1,p-1,q,\dots,k]$. If $p=1$ and $q\leq k-2$, we have $a=W(b)$ for
$b=[\underbrace{1,k-1,\dots,1,k-1}_{2m},1+q,\dots,k]$. This will
complement the existing matching so that the only remaining nontrivial
critical simplices are $[\underbrace{1,k-1,1,\dots,k-1}_{2m},1,k]$, if
$t=km+1$, and $[\underbrace{1,k-1,\dots,1,k-1}_{2m},k]$, if~$t=km$.

It only remains to check that $W$ is still acyclic. Since the newly
matched simplices form an upper ideal of $D_\la$, it is enough to
check the acyclicity condition involving only them. Let
$a_0,\dots,a_f\in D_\la$ be such that $a_0=a_f$, $a_i\neq a_{i+1}$,
$a_i=[\underbrace{1,k-1,\dots,1,k-1}_{2m_i},p_i,q_i,\dots,k]$,
$k-1\geq p_i\geq 2$, and
$$W(a_i)=[\underbrace{1,k-1,\dots,1,k-1}_{2m_i},1,p_i-1,q_i,\dots,k]
\succ a_{i+1},$$ 
for $i=0,\dots,f-1$.

Then, by what we just said, $a_{i+1}$ is obtained from $W(a_i)$ by
merging the blocks indexed $j$ and $j+1$, for $j\geq 2m_i+2$. If
$j\geq 2m_i+3$, or $j=2m_i+2$, but $p_i-1+q_i\neq k-1$, we get into a
contradiction with the choice of~$a_i$'s. Hence we must have
$j=2m_i+2$ and $p_i-1+q_i=k-1$, which implies that $m_i<m_{i+1}$, and
we get a~contradiction $m_0<m_1<\dots<m_f=m_0$.

\noindent
{\bf (b)} The case $\la=(k,2,1^t)$ is very similar. The only
difference is that there is an additional nontrivial critical
$t$-simplex $[2,\underbrace{1,\dots,1}_t,k]$ (according to the idea of
the previous matching, one would want to break $2$ into $1,1$, which
is impossible). Thus, from the previous argument we derive the
conclusion, unless the nontrivial critical cells (for the case
$t+2=km$ and $t+1=km$) are in the neighboring dimensions. These cases
are $(4,2,1,1)$, $(4,2,1,1,1)$, $(3,2,1^4)$, and $(3,2,1^5)$; they can
be verified directly. Observe, that since for some cases we obtain
a homotopy equivalence of $\da_\la$ with a~CW complex with 2 cells
in dimensions higher than 0, we cannot in general determine the
homotopy type of $\da_\la$.
\qed

  Note that since we are not using \cite[Theorem 8.2, Corollary 8.4]{BW} 
for our proof of the Proposition~\ref{pr3.4}(a), we obtain the alternative 
proof of these results of Bj\"orner and Wachs on the homotopy type of the
lattice of intervals generated by all $(k-1)$-element subsets of $\{1,\dots,n-1\}$.

\subsection{Identifications caused by additional resonances}$\,$\label{ss3.2}
\vskip3pt

Let us introduce one more piece of terminology. Recall that, given
a~composition, its resonances are simply linear dependencies of its
parts with coefficients $\pm 1,0$, which we viewed as a~pair of the
subsets of the index set: those parts which get a~coefficient 1 and
those which get a~coefficient -1. We say that a~resonance $r$ can be
{\it derived} from a~set of other resonances $S$ if it follows from
them as a~linear equation. That is, if the entries in an~integer
vector satisfy the equations from the set $S$, then they also satisfy
the equation $r$. A~standard linear-algebraic way to see it is to pass
to the dual vector space and view resonances there as vectors with
$\pm 1,0$ coordinates. Then $r$ can be derived from $S$ if and only if
it lies on the linear span of vectors from $S$. If the resonance $r$
cannot be derived from the set of resonances $S$, we say that $r$
is {\it independent} from $S$.

Having this picture in mind, one can talk about adding independent
resonances to the already existing set. If one could in general
describe what happens to the topology of the corresponding stratum,
then one would have a~general algorithm to compute the algebraic
invariant of $\hyph$. Unfortunately, the combinatorics of the 
situation seems prohibitively complex.

In this subsection we make a~small step on this road.
Let us start with a~simple observation.
\begin{prop}\label{allrel}
  Let $\al=[\al_1,\dots,\al_t]$ and
  $\ti\al=[\ti\al_1,\dots,\ti\al_t]$, $t\geq 3$, be two compositions
  such that the set of resonances of $\ti\al$ is equal to the union of
  the set of resonances of $\al$ with a~new resonance
  $\ti\al_1+\dots+\ti\al_k=\ti\al_{k+1}+\dots+\ti\al_t$, for some
  $1\leq k\leq t-1$.
Let $\la=(\al_1,\dots,\al_t)$ and $\ti\la=(\ti\al_1,\dots,\ti\al_t)$,
then $\wti H_*(\da_{\ti\la})=\wti H_*(\da_\la)\oplus\bz_{(1)}$,
correspondingly $\wti H_*(\widehat{\text{Hyp}_{\ti\la}^n})=
\wti H_*(\hyph)\oplus\bz_{(3)}$.
\end{prop}

\nin{\bf Note}. The condition on the sets of resonances of the
compositions $\al$ and $\ti\al$ in the formulation of the
Theorem~\ref{allrel} is much stronger that just requiring that the
resonance $\sum_{i=1}^k\ti\al_i=\sum_{j=k+1}^t\ti\al_j$ is independent
from the resonances of $\alpha$. It means that this is the {\it only}
resonance added. For example, it implies the following: $\alpha$ has
no resonance of the type $\sum_{i\in I}\al_i= \sum_{j\in J}\al_j$, for
nonempty sets $I$ and $J$, such that simultaneously
$I\subseteq\{1,\dots,k\}$ and $J\subseteq\{k+1,\dots,t\}$.

\vskip3pt \pr Clearly, $\da_{\ti\la}$ is obtained from $\da_\la$ by
gluing together two vertices indexed by the compositions
$[\al_1+\dots+\al_k,\al_{k+1}+\dots+\al_t]$ and
$[\al_{k+1}+\dots+\al_t,\al_1+\dots+\al_k]$. Since the topological
space $\da_\la$ is connected, when $l(\la)\geq 3$, the result follows 
from the homology long exact sequence of a~pair.  
\qed \vskip3pt

  Next we consider the case when the added resonance does not include all
the blocks, but is still the only resonance added.

\begin{sthm2} \label{indrel}
  
Let $[\al_1,\dots,\al_p,\beta_1,\dots,\beta_q,\gamma_1,\dots,\gamma_r]$
be a~composition of type $\lambda$, such that $p\geq 1$, $q\geq 2$,
$r\geq 1$, whose set of resonances includes:
\begin{equation}\label{r1}
\beta_1+\dots+\beta_q=\gamma_1+\dots+\gamma_r.
\end{equation}

Assume there exists $[\ti\al_1,\dots,\ti\al_p,\ti\beta_1,\dots,\ti\beta_q,
\ti\gamma_1,\dots,\ti\gamma_r]$ a~composition of type $\ti\lambda$
such that it has exactly the same resonances as $\la$, except 
for~\eqref{r1} (in particular,~\eqref{r1} is independent from the other
resonances of $\la$). Assume also that the block
$\gamma=\beta_1+\dots+\beta_q+\gamma_1+\dots+\gamma_r$ 
is strongly independent in $(\al_1,\dots,\al_p,\gamma)$.  

Then there exists a~long exact sequence
$$\dots\lra\wti H_{i-2}(\da_{(\al_1,\dots,\al_p)})\stackrel{d_i}{\lra}
\wti H_i(\da_{\ti\la})\lra\wti H_i(\da_\la)\lra
\wti H_{i-3}(\da_{(\al_1,\dots,\al_p)})\stackrel{d_{i-1}}{\lra}\dots\,\,.$$
If, furthermore $\ti\beta_1$ and $\ti\beta_2$ are independent in $\ti\la$, 
then $d_i=0$, and hence $\wti H_i(\da_\la)$ can be found by solving 
the corresponding extension problem. In particular, if 
$\wti H_{i-3}(\da_{(\al_1,\dots,\al_p)})$ is free then
$\wti H_i(\da_\la)=
\wti H_i(\da_{\ti\la})\oplus\wti H_{i-3}(\da_{(\al_1,\dots,\al_p)})$.

\end{sthm2}

\nin{\bf Note 1.} The corresponding information about $\hyph$ can be
derived via the formula $\hyph\cong\susp^2\da_\la$.

\nin{\bf Note 2.} The case $q\geq 1$, $r\geq 2$ is symmetric to the
case considered in the Theorem~\ref{indrel}, hence the same conclusion
can be reached with $\beta$'s and $\gamma$'s interchanged. If $q=r=1$,
then the simplicial complex $\da_\la$ is contractible by the
Stabilization Theorem~\ref{stab}. \vskip3pt

\pr Let $A$ (resp. $\wti A$) be the simplicial subcomplex of
$\da_{\ti\la}$ consisting of the simplices which are labeled by those
compositions, where the sets of blocks
$\{\ti\beta_1,\dots,\ti\beta_q\}$ and
$\{\ti\gamma_1,\dots,\ti\gamma_r\}$ of $\ti\la$ are summed up and the
sum of $\ti\beta_1,\dots,\ti\beta_q$ is either in the same block as
the sum of $\ti\gamma_1,\dots,\ti\gamma_r$ or to the left
(resp.~right) of it. Clearly, $B=A\cap\wti A$ is the simplicial
subcomplex of $\da_{\ti\la}$, where all the blocks
$\ti\beta_1,\dots,\ti\beta_q,\ti\gamma_1,\dots,\ti\gamma_r$ are
summed~up, and $A\cup\wti A$ is the simplicial
subcomplex of $\da_{\ti\la}$, where the sets of blocks
$\{\ti\beta_1,\dots,\ti\beta_q\}$ and
$\{\ti\gamma_1,\dots,\ti\gamma_r\}$ of $\ti\la$ are summed up.

There is a~simplicial map $\da_{\ti\la}\lra\da_\la$ which corresponds
to imposing a~new resonance $\sum_{i=1}^q\ti\beta_i=
\sum_{j=1}^r\ti\gamma_j$ on $\ti\la$. Topologically it corresponds to
gluing the subcomplexes $A$ and $\wti A$ together in the natural way.
There is a~simplicial bijection $\phi$ between $A$ and $\wti A$, which 
interchanges the sums $\sum_{i=1}^q\ti\beta_i$ and
$\sum_{j=1}^r\ti\gamma_j$. This bijection fixes $B$ and therefore we
can glue $A$ together with $\wti A$ inside $\da_{\ti\la}$ by
pointwise identifying those simplices which are mapped to each
other by $\phi$. Gluing together two points in the proof of the 
Proposition~\ref{allrel} is a~simple special case of this procedure.

Let $\bar A$ denote the simplicial subcomplex of $\da_\la$ consisting
of simplices indexed by compositions where the sets
$\{\beta_1,\dots,\beta_q\}$ and $\{\gamma_1,\dots,\gamma_r\}$ are
summed up. We have $\bar A=\da_{(\al_1,\dots,\al_p,\beta,\beta)}$,
where $\beta=\beta_1+\dots+\beta_q=\gamma_1+\dots+\gamma_r$.
The block $\beta$ must be independent in $(\al_1,\dots,\al_p,\beta,\beta)$.

Indeed, since $\gamma$ is strongly independent in 
$(\al_1,\dots,\al_p,\gamma)$, there is no resonance of the type 
$\sum_{i\in I}\al_i+\beta+\beta=\sum_{j\in J}\al_j$. The only
other option for $\beta$ {\it not} to be independent would be to have
a~resonance of the type $\sum_{i\in I}\al_i+\beta=\sum_{j\in J}\al_j$.
But, then $\lambda$ would have resonances 
$\sum_{i\in I}\al_i+\sum_{i=1}^q\beta_i=\sum_{j\in J}\al_j$ and
$\sum_{i\in I}\al_i+\sum_{j=1}^r\gamma_j=\sum_{j\in J}\al_j$.
These two resonances imply~\eqref{r1}, which contradicts to the
existence of $\ti\la$. 

So $\beta$ is independent, and hence, by the Stabilization
Theorem~\ref{stab}(a), $\bar A$ is contractible. Furthermore, it is
clear that $A$, $\wti A$ and $\bar A$ are all isomorphic, hence they
are all contractible.

  By the very nature of the gluing map $\da_{\ti\la}\lra\da_\la$ we have
a~simplicial isomorphism of pairs $(\da_{\ti\la},A\cup\wti A)\cong
(\da_\la,\bar A)$, for the general criteria see Proposition~\ref{pairs}. 
Combining this observation with a~long exact sequence
$$\dots\lra\wti H_i(\bar A)\lra\wti H_i(\da_\la)\lra\wti H_i(\da_\la,\bar A)
\lra\wti H_{i-1}(\bar A)\lra\dots\,\,,$$
we conclude that $\wti H_i(\da_\la)=\wti H_i(\da_\la,\bar A)=
\wti H_i(\da_{\ti\la},A\cup\wti A)$. We also have a~long exact sequence
$$\dots\lra\wti H_i(A\cup\wti A)\stackrel{d_i}{\lra}\wti H_i(\da_{\ti\la})\lra
\wti H_i(\da_{\ti\la},A\cup\wti A)\lra\wti H_{i-1}(A\cup\wti A)
\stackrel{d_{i-1}}{\lra}\dots\,\,.$$

Since both $A$ and $\wti A$ are contractible, we have $\wti H_i(A\cup\wti A)=
\wti H_{i-1}(A\cap\wti A)=\wti H_{i-1}(B)$. Clearly 
$B=\da_{(\al_1,\dots,\al_p,\gamma)}$, and, since $\gamma$ is strongly 
independent in $(\al_1,\dots,\al_p,\gamma)$, we know that 
$B\simeq\susp(\da_{(\al_1,\dots,\al_p)})$,
by the Stabilization Theorem~\ref{stab}(b).

   Let us now see that the homology map $d_i:\wti H_i(A\cup\wti A)\lra
\wti H_i(\da_{\ti\la})$, induced by the inclusion map, is trivial, under 
the condition that $\ti\beta_1$ and $\ti\beta_2$ are independent in $\ti\la$.
Let $K$ be the simplicial subcomplex of $\da_{\ti\la}$ consisting of
the simplices indexed by those compositions, where $\ti\beta_1$ is either in
the same block as $\ti\beta_2$ or in the block with a~smaller index than
the block containing $\ti\beta_2$. 

Clearly, $K\cong\da_{(\bar\al_1,\dots,\bar\al_p,\bar\beta_1,\dots,
  \bar\beta_q,\bar\gamma_1,\dots,\bar\gamma_r)}$, where the set of
resonances of the number partition
$\bar\lambda=(\bar\al_1,\dots,\bar\al_p,\bar\beta_1,\dots,
\bar\beta_q,\bar\gamma_1,\dots,\bar\gamma_r)$ is obtained from the set
of resonances of $\ti\la$ by adding the resonance
$\bar\beta_1=\bar\beta_2$ and everything which it implies together
with the already existing resonances. 

Since $\ti\beta_1$ and $\ti\beta_2$ are independent in $\ti\la$,
$\bar\beta_1=\bar\beta_2$ is independent in $\bar\la$, hence, by the
Stabilization Theorem~\ref{stab}(a), $K$ is contractible. On the other
hand, $K\supseteq A\cup\wti A$, so the inclusion map can be factored
$A\cup\wti A\stackrel{i^1}{\ra} K\stackrel{i^2}{\ra}\da_{\ti\la}$ and
hence $d_*$ factors as well $\wti H_*(A\cup\wti
A)\stackrel{i^1_*}{\ra} \wti H_*(K)\stackrel{i^2_*}{\ra}\wti
H_*(\da_{\ti\la})$. Since the middle term is 0, we conclude that
$d_i=0$. The last conclusion of the theorem now follows.
\qed

\subsection{Noncanonicity of maps between resonances} $\,$
\label{ss3.3}
\vskip3pt

Recall that if $\res\lambda=\res\mu$, then
$\delta_\lambda\simeq\delta_\mu$. Therefore, the simplicial complex
$\delta_{\,\res\lambda}$ is well-defined.

  Whenever we have number partitions $\la>\mu$, there is an~inclusion
map $i(\mu,\la):\da_\mu\ra\da_\la$ which induces $i(\mu,\la)_*:
\wti H_*(\da_\mu)\ra\wti H_*(\da_\la)$, and hence obviously 
 $i(\mu,\la)_*:\wti H_*(\delta_{\,\res\mu})\ra\wti H_*(\delta_{\,\res\lambda})$.
We conjecture however that $i(\mu,\la)_*$ depends on more than just 
the sets $\res\lambda$ and $\res\mu$, i.e., one cannot define 
a unique map $i(\res\mu,\res\la)_*$. More precisely:
\begin{conj}\label{cincl}
  There exist number partitions $\la,\mu,\ti\mu$, such that
\begin{enumerate}
\item $\la>\mu$, $\la>\ti\mu$;
\item $\res\mu=\res\ti\mu$;
\item the homomorphisms of homology groups, $i(\mu,\la)_*$ and $i(\ti\mu,\la)_*$,
induced by the respective inclusion maps are nonisomorphic.
\end{enumerate}
\end{conj}

For every $n\geq 1$, we define a~partial order on the set
$\{\res\lambda\,|\,l(\lambda)=n\}$ by saying that
$\res\ti\la\geq\res\la$ if and only if for some compositions $\alpha$,
resp. $\tilde\alpha$, of type $\lambda$, resp. $\tilde\lambda$, the
set of resonances of $\ti\al$ is a~subset of the set of resonances
of~$\al$. In such a~case, a~choice of the compositions $\alpha$ and
$\tilde\alpha$ induces a~map
$\gamma(\ti\al,\al):\da_{\ti\la}\ra\da_\la$, and further
$\gamma(\ti\al,\al)_*:\wti H_*(\da_{\ti\la})\ra \wti
H_*(\da_\la)$.

\begin{conj}\label{cres}
  The isomorphism type of the map $\gamma(\ti\al,\al)_*$ depends not
  only on the actual number partitions $\ti\la$ and $\la$, rather than
  their sets of resonances, but even on the choice of the pair of
  compositions $\alpha$ and $\tilde\alpha$.
\end{conj}

\nin {\bf Note.} It is easy to come up with examples of number
partitions $\ti\la$ and $\la$, for which there are such pairs of
compositions, which are nonisomorphic under the group of symmetries of
the blocks of $\la$. For example, $\ti\la=(a,b,c,d,d)$,
$\la=(x,x,x,y,y)$, $\tilde\alpha_1=\tilde\alpha_2=[a,b,c,d,d]$,
$\alpha_1=[x,x,x,y,y]$, and $\alpha_2=[x,y,y,x,x]$; the pairs
$(\tilde\alpha_1,\alpha_1)$ and $(\tilde\alpha_2,\alpha_2)$ are
nonisomorphic.

We would like to emphasize that in order to obtain an algorithm for
computing the homology groups of the simplicial complexes
$\delta_{\,\res\lambda}$ it is almost certainly essential to
understand the maps $i(\mu,\la)_*$ and $\gamma(\ti\alpha,\alpha)_*$,
which, as Conjectures~\ref{cincl} and \ref{cres} seem to suggest, may
be a~rather nontrivial task.

\section{Remaining questions and future perspectives}

  We think that understanding the maps $i(\mu,\la)_*$ and 
$\gamma(\phi(\ti\la,\la))_*$ described in Subsection~\ref{ss3.3}, 
combined with the type of the arguments used in the 
proof of the Theorem~\ref{indrel}, would lead to further progress in the 
computation of the homology groups of the simplicial complexes $\da_\la$.

The following observation (proof is left to the reader) is of use when one 
wants to compare two long exact sequences of a~pair, as it was done in 
Theorem~\ref{indrel}.

\begin{prop}\label{pairs}
  Let $(\pi_1,\dots,\pi_{l(\la)})=\la>\mu$ and $(\ti\pi_1,\dots,
\ti\pi_{l(\ti\la)})=\ti\la>\ti\mu$ be number partitions, such that
$l(\la)=l(\ti\la)$ and $l(\mu)=l(\ti\mu)$. We have a~simplicial isomorphism 
of pairs $(\da_\la,\da_\mu)$ and $(\da_{\ti\la},\da_{\ti\mu})$ induced by
$\pi_i\ra\ti\pi_i$, for $1\leq i\leq l(\la)$, (in which case of course 
$\wti H_*(\da_\la,\da_\mu))=\wti H_*(\da_{\ti\la},\da_{\ti\mu})$), if and only 
if the following conditions are satisfied:
\begin{enumerate}
\item $\da_{\ti\mu}$ is the image of $\da_\mu$ under the map induced by 
$\pi_i\ra\ti\pi_i$; 
\item if there is a~resonance $\sum_{i\in I}\pi_i=\sum_{j\in J}\pi_j$ in $\la$, 
then, either we have $\mu\geq(\sum_{i\in I}\pi_i,\sum_{j\in J}\pi_j,\pi_{f_1},
\dots,\pi_{f_t})$, where $\{f_1,\dots,f_t\}=\{1,\dots,l(\la)\}\sm(I\cup J)$, or 
there is a~resonance $\sum_{i\in I}\ti\pi_i=\sum_{j\in J}\ti\pi_j$ 
in $\ti\la$;
\item the same as 2. above, with $\la$ and $\ti\la$, as well as, $\mu$ and $\ti\mu$,
interchanged. 
\end{enumerate}
\end{prop}

{\bf Example of a~computation.}  Let $\la=(3,2,2,1)$, $\mu=(3,3,2)$,
$\ti\la=(5,3,3,1)$, and $\ti\mu=(5,4,3)$, the conditions of the
Proposition~\ref{pairs} are satisfied.  Clearly, $\da_{(3,3,2)}$ is
contractible. By the Stabilization Theorem~\ref{stab}(a) and the
Proposition~\ref{pr3.4} we have $\wti H_*(\da_{(5,3,3,1)})=\bz_{(1)}$,
and, by a~direct observation, $\da_{(5,4,3)}$ is homeomorphic to
$S^1$. Therefore we conclude that $$\wti H_*(\da_{(3,2,2,1)})=\wti
H_*(\da_{(3,2,2,1)},\da_{(3,3,2)})= \wti
H_*(\da_{(5,3,3,1)},\da_{(5,4,3)})=\bz_{(2)}\oplus\bz_{(1)},$$
where
the last equality follows from the fact the the circle $\da_{(5,4,3)}$
does not pass through the vertex indexed by the composition $(6,6)$,
hence the homomorphism of the homology groups $i_*:\wti
H_1(\da_{(5,4,3)})\ra\wti H_1(\da_{(5,3,3,1)})$, induced by inclusion,
is a~zero map.

\end{document}